\documentstyle[12pt]{article}

\newcommand{\sect}[1]{\setcounter{equation}{0}\section{#1}}
\newcommand{\subsect}[1]{\subsection{#1}}

\def\be{\begin{equation}}
\def\ee{\end{equation}}
\def\bea{\begin{eqnarray}}
\def\eea{\end{eqnarray}}

\def\ap{A_+}
\def\am{A_-}
\def\bp{B_+}
\def\bm{B_-}
\def\aa{N}
\def\bb{M}

\def\aap{a_+}
\def\aam{a_-}

\def\ctea{a}
\def\cteb{m}
\def\ctec{b}

\def\discr{{\cal D}_t}

\def\jj{J_3}
\def\jp{J_+}
\def\jm{J_-}

\def\>#1{{\bf #1}}

\def\1{\'{\i}}
\def\R{{\rm I\kern-.2em R}}

\def\sch{{{\cal S}(1+1)}}

\begin{document}

\thispagestyle{empty}

 \
\hfill\

\ 
\vspace{2cm}

\begin{center}
{\LARGE{\bf{Quantum two-photon algebra}}}

{\LARGE{\bf{from non-standard $U_z (sl(2,\R))$ and}}}

{\LARGE{\bf{a discrete time Schr\"odinger equation}}} 
\end{center}

\bigskip\bigskip

\begin{center} Angel Ballesteros$^\dagger$, Francisco J.
Herranz$^\dagger$ and Preeti Parashar$^\ddagger$
\end{center}

\begin{center} {\it {  $^\dagger$ Departamento de F\1sica, Universidad
de Burgos} \\   Pza. Misael Ba\~nuelos, 
E-09001-Burgos, Spain}
\end{center}

\begin{center} {\it {  $^\ddagger$ Department of Mathematics, Texas
A \& M University} \\ College Station, Texas 77843, 
USA}
\end{center}

\bigskip\bigskip\bigskip

\begin{abstract}
The non-standard quantum deformation of the (trivially) extended
$sl(2,\R)$ algebra is used to construct a new quantum deformation of
the two-photon algebra $h_6$ and its associated quantum universal 
$R$-matrix. A deformed one-boson representation for this algebra is
deduced and applied to construct a first order deformation of the 
differential equation that generates the two-photon
algebra eigenstates in Quantum Optics. On the other hand, the
isomorphism between $h_6$ and the (1+1) Schr\"odinger algebra leads to
a new quantum deformation for the latter for which  a
differential-difference realization is presented. From it, a time
discretization of the heat-Schr\"odinger equation is obtained and the
quantum Schr\"odinger generators are shown to be symmetry operators.
\end{abstract}

\bigskip\bigskip\bigskip

\newpage

%%%%%%%%%%%%%%%%% Introduction %%%%%%%%%%%

\sect{Introduction}

The two-photon Lie algebra $h_6$ \cite{Gil} is generated by the
operators
 $\{\aa,\ap,\am,\bp,\bm,\break \bb\}$ endowed with the following
commutation rules
\be
\begin{array}{lll}
 [\aa,\ap]=\ap \quad  &[\aa,\am]=-\am \quad &[\am,\ap]=\bb \cr
[\aa,\bp]=2\bp \quad  &[\aa,\bm]=-2\bm \quad &[\bm,\bp]=4\aa+2\bb
\cr
 [\ap,\bm]=-2\am \quad  &[\ap,\bp]=0 \quad
&[\bb,\,\cdot\,]=0\cr
 [\am,\bp]=2\ap \quad  &[\am,\bm]=0 .\quad
& \cr
\end{array}
\label{aa}
\ee
In this basis, two remarkable Lie subalgebras of $h_6$ can be easily
detected: the harmonic oscillator algebra $h_4$ with generators 
$\{\aa,\ap,\am,\bb\}$,  and the (trivially) extended ${sl}(2,\R)$
algebra generated by $\{\aa,\bp,\bm,\bb\}$ and denoted as
$\overline{sl}(2,\R)$ (the central extension $\bb$ is trivial for  the
latter since it can be absorbed by redefining $\aa\to\aa+\bb/2$).

The two-photon algebra can be used to generate a large zoo of squeezed 
and coherent states for (single mode) one and two-photon processes
which have been recently analysed in \cite{Brif} (in  this respect,
see also \cite{Gil}). More explicitly, the following one-boson
representation of $h_6$
\be
\begin{array}{lll}
\aa=\aap\aam&\qquad \ap=\aap&\qquad \am=\aam\cr
\bb=1&\qquad \bp=\aap^2&\qquad \bm=\aam^2 
\end{array}
\label{ab}
\ee
where
\be
[\aam,\aap]=1 ,
\label{ac}
\ee
shows that one-photon processes are algebraically encoded within the
subalgebra $h_4$ and that $\overline{sl}(2,\R)$ contains the
information concerning two-photon dynamics.
The realization (\ref{ac}) can be translated into a Fock--Bargmann
representation \cite{FB} where the $h_6$ generators act in the
Hilbert space of entire   analytic functions
$f(\alpha)$ as linear differential operators:
\bea
&&\aa=\alpha \frac{d}{d\alpha} \qquad \ap=\alpha \qquad 
\am=\frac{d}{d\alpha}\cr &&\bb=1 \qquad \bp=\alpha^2 \qquad
\bm=\frac{d^2}{d\alpha^2} .
\label{ad}
\eea
The two-photon algebra eigenstates \cite{Brif} are given by the
analytic eigenfunctions that fulfill
\be
(\beta_1\aa + \beta_2 \bm + \beta_3\bp+\beta_4\am+\beta_5\ap)
f(\alpha)=\lambda f(\alpha).
\label{ae}
\ee
In the Fock--Bargmann
representation (\ref{ad}), the
following differential equation is deduced from (\ref{ae}):
\be
\beta_2 \frac{d^2f}{d\alpha^2}+(\beta_1\alpha  +\beta_4)\frac{d
f}{d\alpha} +(\beta_3\alpha^2+\beta_5\alpha -\lambda)f=0 
\label{af}
\ee
where $\beta_i$ are arbitrary complex coefficients and  $\lambda$ is a
complex eigenvalue. The solutions of this equation  (provided a
suitable normalization is imposed) give rise to the two-photon
coherent/squeezed states \cite{Brif}. One and two-photon coherent and
squeezed states corresponding to the subalgebras $h_4$ and
$\overline{sl}(2,\R)$ can be derived from equation (\ref{af}) by
setting $\beta_2=\beta_3=0$ and  $\beta_4=\beta_5=0$, respectively.

A novel (to our knowledge) and interesting feature of $h_6$ has been
recently pointed out in \cite{tufoton}: if we define 
\be
D=-\aa-\frac 12\bb \qquad P=\ap\qquad K=\am\qquad
H=\frac 12\bp\qquad C=\frac 12\bm ,
\label{fb}
\ee
and compute the commutation rules in this new basis we obtain 
\be
\begin{array}{lll}
 [D,P]=-P \quad  &[D,K]=K \quad &[K,P]=M \cr
 [D,H]=-2H \quad  &[D,C]=2C  \quad &[H,C]=D \cr
 [K,H]=P \quad  &[K,C]=0 \quad
&[M,\,\cdot\,]=0\cr
 [P,C]=-K \quad  &[P,H]=0 ,\quad
& \cr
\end{array}
\label{fa}
\ee
which are just the defining relations of the (1+1) dimensional
Schr\"odinger  algebra ${\cal S}(1+1)$ \cite{Hagen,Nied} generated
by the time translation operator $H$, space translation $P$, Galilean
boost $K$, dilation $D$, conformal transformation $C$  and mass $M$.
Therefore, the
isomorphism (\ref{fb}) relates in a simple way two physically different
frameworks and enhances the role played by the underlying abstract Lie
symmetry. The kinematical nature of ${\cal S}(1+1)$ can be
explicitly put forward by recalling the usual differential realization
of the Schr\"odinger generators in terms of the time and space
coordinates $(t,x)$:
\bea
&&H=\partial_t \qquad  P=\partial_x \qquad M=\cteb \cr
&&K=-t \partial_x - m x\qquad D=2 t \partial_t +
x\partial_x -\ctea \label{ffa}\\
&&C=t^2 \partial_t + t x \partial_x - \ctea t  + \cteb x^2/2 
\nonumber
\eea
where $\ctea$, $\cteb$ are the constants that label the
representation. This algebra appears as the symmetry algebra for the
(1+1) heat/Schr\"odinger equation (SE), that can be straightforwardly
obtained from the representation (\ref{ffa}) as the Casimir operator
$P^2 -2\,M\,H$ of the (1+1) Galilei subalgebra of ${\cal S}(1+1)$.

The direct connection of ${\cal S}(1+1)$ with the SE has
recently motivated the search for $q$-deformations of this algebra
connected with discretizations (in both space and time variables) of
this equation on geometric lattices of the type $y_n=q^n\,y_0$. In
\cite{Vinet,FV} $q$-deformations of the vector field realization
(\ref{ffa}) arise as symmetry algebras of different discretized
versions of the SE. From a constructive point of view, another
$q$-Schr\"odinger algebra has been directly introduced in \cite{Dobrev}
in order to get a generalized $q$-SE from its deformed representation
theory. However, none of these $q$-algebras has been found to be
consistent with a Hopf algebra structure \cite{CP} yet. On the other
hand, a discretization of the SE on a regular space-time lattice 
$y_n=y_0+n\,z$ has been also studied from a symmetry approach in 
\cite{javier}. In this case, the resultant symmetry operators close a
non-deformed Schr\"odinger algebra. This classical nature of the
discrete symmetries on a regular grid has been also proven for the
wave equation in \cite{JavLuis}. 

The first Hopf algebra 
deformation of $h_6$ 
and the isomorphism (\ref{fb}) used to derive the
corresponding quantum Schr\"odinger algebra relations have been
recently introduced in \cite{tufoton}. This quantum algebra was
obtained by starting from the (non-standard) quantum deformation of
its one-photon subalgebra
$h_4$ \cite{osc}.  This paper is devoted to the
investigation of a new quantum deformation of $h_6\equiv {\cal
S}(1+1)$ coming from the (non-standard) quantum deformation of
the other relevant subalgebra $\overline{sl}(2,\R)$ given in
\cite{bosones}. As a result, the corresponding deformed analogues of
(\ref{ad}) and (\ref{ffa}) are shown to provide two systematic
applications of this deformed symmetry: among them, we shall
insist in the construction of a differential-difference analogue of
the SE with quantum algebra symmetry. 

In the next
section the explicit form of the quantum two-photon algebra is
presented starting from the classical $r$-matrix $r=z \aa\wedge \bp$ 
that underlies the deformation. Some mathematical properties are
stated: among them, the embedding $U_z(\overline{sl}(2,\R))\subset
U_z(h_6)$ and the proof that the universal $R$-matrix of
$U_z(\overline{sl}(2,\R))$ is also valid for $U_z(h_6)$.  A deformed
one-boson representation is afterwards presented and
translated into a Fock--Bargmann realization. This representation
suggests the construction of deformed states of light through a
differential equation (whose classical version is (\ref{af})) for
their corresponding quantum optical analytic functions. Since the
general equation so obtained is rather involved, a truncated
deformation (valid for small values of $z$) is explicitly presented. 

The third section is devoted to the analysis of the quantum
Schr\"odinger algebra originated from the previous deformation through
the isomorphism (\ref{fb}). A systematic construction of a deformed
SE is presented from the quantum deformation of the
representation (\ref{ffa}), that is obtained in terms of 
differential-difference operators. The fact that the Galilei
generators $\{K,H,P,M\}$ close a (deformed) subalgebra allows us to
consider the corresponding Casimir operator, that leads to a  time
discretization of the SE on a uniform time
lattice. Finally, the full quantum algebra invariance of the equation
is proven by checking that the remaining generators $D$ and $C$ also
transforms solutions into solutions. This result can be interpreted
as a first example of uniform-lattice discretization induced from
quantum algebras. From a physical point of view, an approach to the
discrete time dynamics on such a time lattice has been recently
developed in \cite{JN}.

%%%%%%%%%%%%%!%%% quantum two-photon algebra %%%%%%%%%%%

\sect{A new quantum two-photon algebra}

\subsect{The quantum algebra $U_z(h_6)$}

A general property of non-standard classical $r$-matrices linked to
a given Lie algebra $g$ is their automatic compatibility with any
other Lie algebra $g'$ that contains $g$ as a subalgebra. Therefore,
some Lie bialgebra structures on $g'$ can be obtained by using
(known) classical $r$-matrices of $g$. If we recall that Lie
bialgebras are just the first order of quantum deformations \cite{CP},
it turns out that the construction of (non-standard) quantum algebra
deformations of
$g'$  can be guided by known results for their subalgebras. 

This was the case in
\cite{tufoton}, where the non-standard
quantum  deformation of $h_4$ introduced in \cite{osc} was used
to obtain a so  called Jordanian quantum two-photon algebra such that
$U_z(h_4)\subset U_z(h_6)$; the classical $r$-matrix
linked to that deformation was $r=z \aa\wedge \ap $, and the primitive
generators in $h_6$ turned out to be $\bb$ and $\ap$. 

The same kind
of construction can be obtained by means of the  classical $r$-matrix
of $sl(2,\R)$, the subalgebra of $h_6$ linked to ``pure" two-photon
processes:
\be
r= z \aa\wedge \bp.
\label{ba}
\ee
By construction, (\ref{ba}) is a solution of the classical Yang--Baxter
equation for $h_6$   and  generates  the Lie bialgebra with
cocommutators, $\delta(X)=[1\otimes X + X\otimes 1,r]$,  given by:
\be
\begin{array}{l}
 \delta(\bp)=0\qquad \delta(\bb)=0\cr
 \delta(\aa)=2z\, \aa\wedge \bp\qquad 
 \delta(\ap)=-z\, \ap\wedge \bp\cr
 \delta(\am)=z\, (\am\wedge \bp + 2 \aa\wedge \ap)\cr
 \delta(\bm)=2z\, (\bm\wedge \bp + \aa\wedge \bb) . 
\end{array}
\label{bb}
\ee
The four generators $\{\aa,\bp,\bm,\bb\}$
close a sub-bialgebra which is isomorphic to the one  
described for $\overline{sl}(2,\R)$ in \cite{bosones} with  the aid
of the  identification: 
\be
\jp = \frac {\bp}2\qquad  \jm=-\frac {\bm}2\qquad  \jj=\aa   \qquad 
I=-\frac {\bb}2  
\label{bc}
\ee
and the replacement $z\to 2z$. Thus we  can profit  from
this fact getting a new   quantum two-photon algebra with different
properties with respect to the Jordanian one  introduced in
\cite{tufoton}.

The steps  involved  in the construction of the quantum deformation of
the above two-photon bialgebra are as follows:
{\it i)} Take  the Hopf structure
of $U_z(\overline{sl}(2,\R))$ written in the two-photon basis. {\it
ii)} Find a coassociative coproduct $\Delta$
of the   remaining two  generators $\ap$ and $\am$ by taking into
account that (\ref{bb}) gives its first order deformation.  {\it iii)}
 Compute the remaining deformed commutation rules by imposing the
coproduct to be a Lie algebra homomorphism.  {\it iv)} Finally, counit
$\epsilon$ and antipode $\gamma$ are deduced from the previous
results.

The resulting  quantum Hopf algebra $U_z(h_6)$ reads:
\be
\begin{array}{l}
\Delta(\bp)=1\otimes \bp + \bp \otimes 1 \qquad
\Delta(\bb)=1\otimes \bb + \bb \otimes 1\cr
\Delta(\aa)=1\otimes \aa + \aa \otimes e^{2z\bp}
\qquad \Delta(\ap)=1\otimes \ap + \ap \otimes
e^{-z\bp} \cr
\Delta(\am)=1\otimes \am + \am \otimes e^{z\bp}
+ 2z\aa \otimes e^{2z\bp} \ap\cr
\Delta(\bm)=1\otimes \bm + \bm \otimes e^{2z\bp}
+ 2z\aa \otimes e^{2z\bp} \bb\cr
\end{array}
\label{bd}
\ee
\be
\epsilon(X) =0 \qquad
 \mbox{for\ \ $X\in \{\aa,\ap,\am,\bp,\bm,\bb\}$}
\label{be} \ee
\be
\begin{array}{l}
\gamma(\ap)=-\ap e^{z\bp}\qquad \gamma(\bb)=-\bb\cr
\gamma(\aa)=-\aa e^{-2z\bp}\qquad
\gamma(\bp)=-\bp \cr
\gamma(\am)=-(\am -2z\aa\ap) e^{-z\bp}\cr
\gamma(\bm)=-(\bm -2z\aa\bb) e^{-2z\bp}\cr
\end{array}
\label{bf}
\ee
\bea
&& [\aa,\ap]=\ap \qquad   [\aa,\am]=-\am \qquad
 [\am,\ap]=\bb   \cr
&& [\aa,\bp]=\frac {e^{2z\bp}-1}z \qquad    [\aa,\bm]=-2\bm - 4z
{\aa}^2 \cr
&&[\bm,\bp]= 4\aa + 2\bb e^{2z\bp} \qquad  [\bb,\,\cdot\,]=0 
\label{bg}\\
 && [\ap,\bm]=-2\am + 2z(\aa \ap + \ap \aa) \qquad
[\ap,\bp]=0\cr
&&  [\am,\bp]=2 e^{2z\bp} \ap
  \qquad   [\am,\bm]=-2z(\aa \am + \am \aa).
\nonumber
\eea

Clearly by construction we have $U_z(\overline{sl}(2,\R))\subset  
U_z(h_6)$.  Note that the  oscillator algebra $h_4$ is  preserved as an
undeformed subalgebra  only at the level of commutation relations.

%%%%%%%%%%%%%%%%%  universal $R$-matrix %%%%%%%%%%%

\subsect{Universal $R$-matrix}

The quantum universal $R$-matrix for $U_z(h_6)$ is given by
 \be
{\cal R}=\exp\{-z\bp\otimes\aa\}\exp\{z\aa\otimes\bp\} ,
\label{ca}
\ee
 which is a solution of the quantum Yang--Baxter equation
and verifies the property
\be
{\cal R}\Delta(X){\cal R}^{-1}=\sigma\circ \Delta(X) \qquad
\mbox{for}\ \ X\in h_6 
\label{cb}
\ee
where $\sigma (a\otimes b)=b\otimes a$.

In order to prove this statement we have to take into account that
  the element (\ref{ca}) is a universal $R$-matrix for
$U_z(\overline{sl}(2,\R))$ \cite{bosones}.  Therefore, as  
$U_z(\overline{sl}(2,\R))\subset U_z(h_6)$ we  only need to show  
that  $\cal R$ also satisfies (\ref{cb}) for the remaining two
generators $\ap$ and $\am$.

We consider the formula
\be
e^{f}\,\Delta(X)\,e^{-f}=\Delta(X) +\sum_{n=1}^\infty \frac
1{n!}\,[f,[\dots [f,\Delta(X)]]^{n)}\dots]  
\label{cj}
\ee
and set   $X\equiv \ap$, $f\equiv z\aa\otimes\bp$.  A direct
computation shows that
\be
[z\aa\otimes\bp,[\dots [z\aa\otimes\bp,\Delta(\ap)]]^{n)}\dots]
=\ap\otimes  e^{-z\bp}(z\bp)^n\quad n\ge 1
\ee
so that
\bea
&&e^{z\aa\otimes\bp}\Delta(\ap)e^{-z\aa\otimes\bp}
=\Delta(\ap)
+\sum_{n=1}^\infty \ap\otimes  e^{-z\bp}\frac {(z\bp)^n}{n!}\cr
 &&\qquad = 1\otimes \ap + \ap \otimes e^{-z\bp}+
\ap \otimes e^{-z\bp}(e^{z\bp}-1)\cr
&&\qquad =1\otimes \ap + \ap \otimes
1\equiv \Delta_0(\ap) .\label{ck}
\eea
On the other hand we find
\be
[-z\bp\otimes\aa,[\dots [-z\bp\otimes\aa,\Delta_0(\ap)]]^{n)}\dots]
=(-z\bp)^n  \otimes \ap \quad n\ge 1
\ee
and the proof for $\ap$ follows
\bea
&&e^{-z\bp\otimes\aa}\Delta_0(\ap)e^{z\bp\otimes\aa}
=\Delta_0(\ap)
+\sum_{n=1}^\infty \frac {(-z\bp)^n}{n!} \otimes \ap\cr
 &&\qquad = 1\otimes \ap + \ap \otimes1 +
(e^{-z\bp}-1)  \otimes \ap =\sigma\circ \Delta(\ap)  .\label{cl}
\eea

Now  let $X\equiv \am$  and  $f\equiv
z\aa\otimes\bp$:
\be
[z\aa\otimes\bp, \Delta(\am)]
=-2z\aa\otimes e^{2z\bp}\ap  - z\am\otimes \bp
e^{z\bp}
\ee
\be
[z\aa\otimes\bp,[\dots [z\aa\otimes\bp,\Delta(\am)]]^{n)}\dots]
=\am\otimes  e^{z\bp}(-z\bp)^n\quad n\ge 2 .
\ee
Then
\bea
&&\!\!\!\! e^{z\aa\otimes\bp}\Delta(\am)e^{-z\aa\otimes\bp}
=\Delta(\am) -2z\aa \otimes e^{2z\bp}\ap
+\sum_{n=1}^\infty \am\otimes  e^{z\bp}\frac {(-z\bp)^n}{n!}\cr
 &&\qquad = 1\otimes \am + \am \otimes
1\equiv \Delta_0(\am) .\label{cp}
\eea
Similarly we find that
\be
[-z\bp\otimes\aa, \Delta(\am)]
= 2z e^{2z\bp}\ap\otimes \aa + z\bp\otimes \am    
\ee
\be
[-z\bp\otimes\aa,[\dots [-z\bp\otimes\aa,\Delta_0(\ap)]]^{n)}\dots]
=(z\bp)^n  \otimes \am \quad n\ge 2 
\ee
and, consequently,
\bea
&&\!\!\!\! e^{-z\bp\otimes\aa}\Delta_0(\am)e^{z\bp\otimes\aa}
=\Delta_0(\am) +2z e^{2z\bp}\ap \otimes \aa
+\sum_{n=1}^\infty \frac {(z\bp)^n}{n!} \otimes \am\cr
 &&\qquad = \am \otimes 1 +
e^{z\bp} \otimes \am +2z e^{2z\bp}\ap \otimes \aa   =\sigma\circ
\Delta(\am) .\label{cs}
 \eea

%%%%%%%%%%%%%%%%%  one-boson representation %%%%%%%%%%%

\subsect{Deformed Fock-Bargmann realization}

 Let  $\aam$, $\aap$  be  the  boson generators fulfilling (\ref{ac}).
Then a one-boson representation of $U_z(h_6)$ (\ref{bg}) is given by: 
\bea
&&\bp=\aap^2\qquad \bb=1\qquad
 \aa=\frac{e^{2z\aap^2}-1}{2z\aap}\,\aam \cr
&& \ap=\left(\frac{1-e^{-2z\aap^2}}{2z}\right)^{1/2}\qquad 
 \am=\frac{e^{2z\aap^2}}{\aap}
\left(\frac{1-e^{-2z\aap^2}}{2z}\right)^{1/2}\aam \label{da}\\
&& \bm =\left( \frac{e^{2 z \aap^2}-1}{2 z \aap^2 }   \right)   
\aam^2 + 
\left(\frac {e^{2 z \aap^2}}{\aap}+\frac {1-e^{2 z \aap^2}} {2 z
\aap^3}\right)  
\aam .
\nonumber
\eea
Note that in the limit $z\to 0$ we recover the classical representation
(\ref{ab}).  Let us also remark that, in spite of the fact that
(\ref{bg}) presents a non-deformed oscillator subalgebra in terms of
the abstract generators $\{\aa,\ap,\am,\bb\}$, their representation
(\ref{da}) includes strong deformations in terms of usual boson
operators. The above  deformed one-boson representation provides  the
Fock--Bargmann one by  setting $\aap\equiv \alpha$ and $\aam\equiv
\frac d{d\alpha}$. However, due to its rather complicated form we will
restrict  to the first order in $z$. In such approximation, the  
two-photon differential operators   acting on the space of entire
analytic functions
$f(\alpha)$ turn out to be:
\bea
&&\bp=\alpha^2\qquad \bb=1\qquad
 \aa=(\alpha+z\alpha^3)\frac{d}{d\alpha}\cr
&&\ap=\alpha-  z \alpha^3/2    \qquad   
\am=(1+3z\alpha^2/2)\frac{d}{d\alpha} \label{db}\\
&&\bm=(1+z\alpha^2)\frac{d^2}{d\alpha^2}+z\alpha\frac{d}{d\alpha} .
\nonumber
\eea
 We substitute these operators in (\ref{ae}) getting the  differential
equation:
\bea
&&\beta_2(1+z\alpha^2)  \frac{d^2f}{d\alpha^2}+\left(\beta_1\alpha
+\beta_4
+z(\beta_1\alpha^3+\beta_2\alpha+3\beta_4\alpha^2/2)\right)\frac{d
f}{d\alpha}\cr &&\qquad +(\beta_3\alpha^2+\beta_5\alpha
-z\beta_5\alpha^3/2-\lambda)f=0 .
\label{dc}
\eea
The particular equation with $\beta_4=\beta_5=0$  is associated to the
Hopf subalgebra $U_z(\overline{sl}(2,\R))$ while the case
$\beta_2=\beta_3=0$ corresponds to the harmonic oscillator sector
(recall that the latter is not a Hopf subalgebra). Besides the
problem of finding particular solutions of the equation (\ref{dc}),
there are two other important questions to be considered:
it is necessary to define a deformed `invariant' measure and  the
possible solutions must be analytic.

%%%%%%%%%%%%%%%%% quantum Schr\"odinger algebra %%%%%%%%%%%

\sect{Quantum Schr\"odinger algebra and time discre\-tization}

As it was mentioned in the Introduction, a precise Lie subalgebra
of ${\cal S}(1+1)$ can  be considered in order to explain the
significance of this algebra as the Lie symmetry algebra of the (1+1) 
SE: the extended (1+1) Galilei algebra generated by
$\{K,H,P,M\}$. In fact, the (free) heat-Schr\"odinger equation  \be
(\partial_x^2 - 2 m \partial_t) \phi(x,t)=0, 
\label{sch}
\ee
is just the Casimir operator $P^2-2 M H$ of the Galilei subalgebra
in the representation (\ref{ffa}). If we denote the equation operator
as
\be
E=P^2-2 M H,
\label{equation}
\ee
we can say that an operator $S$ is a symmetry of $E$ if $S$
transforms solutions into solutions:
\be
(E\,S) \phi(x,t)= (\Lambda\,E)\phi(x,t),
\label{symm}
\ee
where $\phi(x,t)$ is a solution of (\ref{sch}) and $\Lambda$ is
another operator. In particular, since $E$ is given by the Casimir
operator of the extended Galilei algebra, all its generators are
symmetries $S$ of the SE fulfilling (\ref{symm})
with $\Lambda=S$. Moreover, the dilation operator $D$ is also a
symmetry, since (\ref{fa}) implies that $[E,D]=2 E$. Finally, we can
compute
\be
[E,C]=-(K P + P K +2 M D).
\label{conf}
\ee
Now, if we take the particular representation with $\ctea=-1/2$, we
obtain that the right hand side of (\ref{conf}) is just the operator
$2 t E$, so that the conformal operator $C$ also
transforms solutions of $E$ into solutions. As we shall see in the
sequel, this procedure can be generalized to the quantum case.

\subsect{Quantum Schr\"odinger algebra}

Let us now use the isomorphism (\ref{fb}) in order to obtain
another  non-standard  quantum deformation of $\sch$. The
new classical $r$-matrix and cocommutators read \be
r=2z H\wedge D +z H\wedge M
\label{fc}
\ee
\be
\begin{array}{l}
\delta(H)=0\qquad \delta(M)=0\qquad
\delta(P)=-2z P\wedge H\cr
 \delta(K)=z(2K\wedge H +2 P\wedge D + P\wedge M)\cr
\delta(D)=z(4D\wedge H + 2 M\wedge H)\qquad
 \delta(C)=z(4C\wedge H- D\wedge M) .
\end{array}
\label{fd}
\ee 
The coproduct and the commutation rules of the Hopf algebra
$U_z(\sch)$ are
\bea
&&\Delta(H)=1\otimes H + H \otimes 1
\qquad \Delta(M)=1\otimes M + M \otimes 1\cr
&&\Delta(P)=1\otimes P + P\otimes e^{-2zH}\cr
&&\Delta(K)=1\otimes K + K\otimes e^{2zH}-2z(D+  M/2)\otimes
e^{4zH}P\cr
&&\Delta(D)=1\otimes D + D\otimes e^{4zH} +   M\otimes
(e^{4zH}-1)/2\cr
&&\Delta(C)=1\otimes C + C\otimes  e^{4zH}-z (D +  
M/2)\otimes  e^{4zH}  M,
\label{fe}
\eea
\bea
&& [D,P]=-P \qquad   [D,K]=K\qquad
 [K,P]=M   \cr
&& [D,H]= \frac {1-e^{4zH}}{2z} \qquad   [D,C]=2 C +2z  (D+ 
M/2)^2\cr 
&&[H,C]=  D+ M (1-e^{4zH})/2\qquad  
  [K,H]= e^{4zH}P\cr
&& [K,C]=z(DK+KD+KM)\qquad [P,H]=0\cr
&&  [P,C]=-K-z  (DP+PD+PM)  \qquad 
 [M,\,\cdot\,]=0 . 
\label{ff}
\eea

Notice that in this deformation   the  primitive
generator  (besides the mass) is the time translation $H$, while
in \cite{tufoton} this role was played by the space translation $P$.
The quantum $\overline{sl}(2,\R)$ algebra corresponds now to the Hopf
subalgebra spanned by $\{D,H,C,M\}$, and the extended Galilei
generators together with the dilation, $\{K,H,P,M,D\}$, also close a
 Hopf subalgebra.

On the other hand, the 
universal $R$-matrix for
$U_z(\sch)$ is provided by (\ref{ca}):
\be
\begin{array}{l}
{\cal R}=\exp\{2zH\otimes D\}\exp\{z H\otimes M \}
\exp\{-z M\otimes H \}\exp\{-2z D\otimes H\} .\cr
\end{array}
\ee

\subsect{A time discretization of the Schr\"odinger equation}

A realization of $U_z(\sch)$ (with classical limit  (\ref{ffa})) 
reads
\bea
&&H=\partial_t \qquad  P=\partial_x \qquad M=\cteb \cr
&&K=- (t + 4z) e^{4 z \partial_t}  \partial_x  - m x  \cr  
&& D=2 (t+ 4z) \frac{e^{4 z \partial_t}-1}{4z} +
x\partial_x -\ctea   \label{fg}\\
&&C=(t^2 - 4z \ctec t) \frac{e^{4 z \partial_t}-1}{4z}  + t x
\partial_x - \ctea t  + \cteb x^2/2    - 4 z (\ctec +1) e^{4 z
\partial_t}  \cr
&&\qquad - z x^2 \partial_{x}^2  -2 z (\ctec-\ctea + 1/2) x
\partial_x -z(\ctec-\ctea)^2
\nonumber
\eea
where $\ctec= m/2-2$.  If we introduce a discrete time derivative as
\be
\discr f(t,x):=\left(\frac{e^{4 z \partial_t}-1}{4z}
\right) f(t,x)=\frac{ f(t+4z,x)-f(t,x)}{4z},
\label{fh}
\ee
 then   (\ref{fg}) adopts the form of a differential-difference
realization:
\bea
&&H=\partial_t \qquad  P=\partial_x \qquad M=\cteb \cr
&&K=- (t + 4z)(1+ 4z \discr)  \partial_x  - m x 
\cr   
&& D=2 (t+ 4z) \discr +
x\partial_x -\ctea   \label{fi}\\
&&C=(t^2 - 4z \ctec t) \discr  + t x
\partial_x - \ctea t  +
\cteb x^2/2     - 4 z (\ctec +1)  (1+ 4z \discr)  \cr
&&\qquad  - z x^2 \partial_{x}^2  -2
z (\ctec-\ctea + 1/2) x
\partial_x -z(\ctec-\ctea)^2
\nonumber
\eea

We can now try to reproduce the symmetry procedure before outlined
for the non-deformed Schr\"odinger algebra. The first important
feature of the deformation just introduced is that, at the level of
commutation rules, the Galilei structure remains as a
(deformed) subalgebra. It is not difficult to check that the
deformed Casimir operator for that subalgebra is
\be
E_z =P^2-2 M \frac {1-e^{-4zH}}{4z}.
\label{defcas}
\ee
It seems quite natural to define the corresponding SE  as the action
of (\ref{defcas}) on $\phi(x,t)$ through (\ref{fg}). It reads
\be
  \left(\partial_x^2 - 2 m \frac
{1-e^{-4z\partial_t}}{4z}\right)
\phi(x,t)=0\quad \equiv \quad  (\partial_x^2 - 2 m {\cal D}_t^{-})
\phi(x,t)=0,
\label{qsch}
\ee
which is just a time discretization of the free SE. The discrete time
derivative ${\cal D}_t^{-}$  is related to (\ref{fh}) by changing
$z\to -z$. 

By construction, the quantum algebra generators $\{K,H,P,M\}$ are
symmetries of $E_z$. Furthermore, $D$ and $C$  are also deformed
Schr\"odinger symmetries. This fact follows straightforwardly for the
dilation since from (\ref{ff})  
\be
[E_z,D]= 2 E_z .
\ee
For the conformal transformation, we have
\bea
&&[E_z,C]=-(KP+PK + 2 M D e^{-4 z H}) + M(M+2) (1-e^{-4 z H})\cr
&&\qquad\qquad\quad - z(DP^2 + 2 PDP + P^2 D+2 P^2 M) .
\eea
By introducing the realization (\ref{fg}) we find that
\bea
&&[E_z,C]=2(t+ 4z - 2 z x \partial_x) E_z - 2 z \left\{ m + 2
(1-a)\right\}\partial_x^2\cr
&&\qquad\qquad\quad +m(1+ 2 a e^{-4 z \partial_t})
+ m (m+2 ) ( 1 - e^{-4 z \partial_t}) .
\eea
Now, as in the non-deformed case, we have to set $a=-1/2$ in order to
obtain $C$ as a symmetry:
\be
[E_z,C]=2\left\{ t+  z (1- m - 2   x \partial_x)\right\} E_z .
\ee

\sect{Concluding remarks}

The connection between $q$-algebras and geometric $q$-lattices seems to
be a general fact \cite{FV}, although a general theory of this
kind of systems from the point of view of symmetries is still
lacking. The results presented in this paper and the structural
properties of non-standard deformations suggest a natural link
between non-standard quantum algebras and uniform lattices whose
generality and physical implications should be better understood. 

In this sense, it is remarkable that the quantum
Schr\"odinger algebra here presented induces in a very simple way a
natural definition of a  discrete time SE that is
invariant under the quantum algebra generators. It
becomes apparent that a similar construction  for the non-standard
deformation given in \cite{tufoton} would lead to a
(uniform) discretization of the SE in the spatial
coordinate. This different discretization would be, in the kinematical
context, the algebraic consequence of the use of either the $sl(2,\R)$
or the $h_4$ subalgebra as generating structures for the construction
of the deformation.

Other discretizations of the  SE linked with other quantum algebras
can be found in 
\cite{ita,DD,Micu}. In particular, the use in \cite{DD} of
the $sl(2,\R)$ symmetry in order to include some potentials was also
analysed in \cite{FV} in a $q$-lattice context and would deserve
further attention from the point of view of our deformation.
Another open problem arising from the results here presented is the
study of the special functions that can be defined from them.
Work in all these directions is in progress.

%\newpage

%%%%%%%%%%%%%%%%% ACKNOWLEDGEMENTS %%%%%%%%%%%

\bigskip
\bigskip

\noindent
{\Large{{\bf Acknowledgments}}}

\bigskip

The authors acknowledge J. Negro for helpful suggestions. A.B. and
F.J.H. have been partially supported by DGICYT (Project  PB94--1115)
from the Ministerio de Educaci\'on y Cultura  de Espa\~na and by Junta
de Castilla y Le\'on (Projects CO1/396 and CO2/297). P.P. thanks
Prof. E. Letzter for hospitality at Texas A \& M University.

\end{document}